  \documentclass[reqno,11pt]{article}
\usepackage{amsfonts,latexsym}
\usepackage{amsmath}
\usepackage{amscd}
\usepackage{float,amsmath,amssymb,mathrsfs,bm,multirow,graphics}
\usepackage[dvips]{graphicx}
\usepackage{todonotes}
\usepackage[percent]{overpic}

\usepackage[
  style=nature,
  sorting = nyt,
  maxalphanames=4]{biblatex}
\newcommand{\R}{{\mathbb R}}

\newcommand{\C}{{\mathbb C}}

\newcommand{\reg}{\mathrm{reg}}

\newcommand{\res}{\text{\upshape Res\,}}
\newcommand{\re}{\text{\upshape Re\,}}
\newcommand{\im}{\text{\upshape Im\,}}
\newcommand{\erf}{\text{\upshape erf\,}}

\newcommand{\proofbegin}{\noindent{\it Proof.\quad}}
\newcommand{\proofend}{\hfill$\Box$\bigskip}

\def\Xint#1{\mathchoice
{\XXint\displaystyle\textstyle{#1}}%
{\XXint\textstyle\scriptstyle{#1}}%
{\XXint\scriptstyle\scriptscriptstyle{#1}}%
{\XXint\scriptscriptstyle\scriptscriptstyle{#1}}%
\!\int}
\def\XXint#1#2#3{{\setbox0=\hbox{$#1{#2#3}{\int}$}
\vcenter{\hbox{$#2#3$}}\kern-.5\wd0}}

\def\dashint{\;\Xint-}

\newtheorem{theorem}{Theorem}[section]

\newtheorem{lemma}[theorem]{Lemma}

\newtheorem{remark}[theorem]{Remark}

\newtheorem{figuretext}{Figure}

\numberwithin{equation}{section}

\usepackage[colorlinks=true]{hyperref}
\hypersetup{urlcolor=blue, citecolor=red, linkcolor=blue}

\input epsf
\title{\sc A novel asymptotic technique for integrals involving the Hankel contour and the Bleistein asymptotic formula}
\author{A. S. Fokas$^{1}$ and J. Lenells$^{2}$}
\date{{\small $^1$Department of Applied Mathematics and Theoretical Physics, University of Cambridge, Cambridge CB3 0WA, United Kingdom
	\\
$^{2}$Department of Mathematics, KTH Royal Institute of Technology, \\ 100 44 Stockholm, Sweden}}

\begin{document}
\maketitle
\begin{abstract} 
\noindent
Several important functions, including the gamma function, as well as several infinite sums, admit integral representations involving the Hankel contour. In addition, the large $t$ asymptotic analysis of several recently derived identities satisfied by the Riemann zeta function requires computing the asymptotic form of certain integrals which also involve the Hankel contour; these integrals depend on a real parameter, $\alpha$. A rigorous asymptotic technique is presented here for computing such integrals to all orders. For certain values of $\alpha$, the relevant formula, in addition to an asymptotic series of explicit terms, also contains a specific integral. It is shown that, remarkably, the leading order of this integral can be written in the form of the leading order of the Bleistein integral. The latter integral arises in the implementation of the classical steepest descent method in the case that the stationary point coincides with one of the boundary points of the integral under consideration.
 \end{abstract}




\section{Introduction}\label{introsec}

\begin{figure}
\begin{center}
\begin{overpic}[width=.5\textwidth]{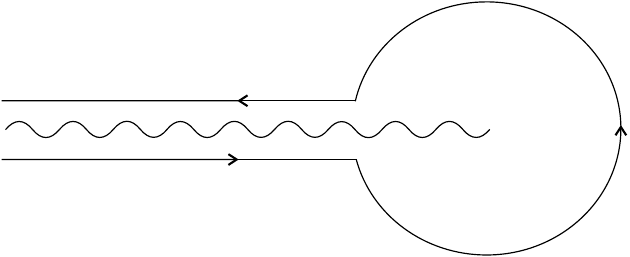}
     \put(101,18){$H_1$}
     \put(80,18){$0$}
\end{overpic}
     \begin{figuretext}\label{Hankel.pdf}
         The Hankel contour $H_1$ surrounds the negative real axis in the counterclockwise direction.
     \end{figuretext}
  \end{center}
\end{figure}

Several important functions and infinite sums admit integral representations which involve the Hankel contour, $H_1$, depicted in Figure \ref{Hankel.pdf} and defined by
\begin{align}\label{}
H_1 = \left\{ re^{-i\pi} : 1<r<\infty \right\} \cup \{ e^{i\theta} : -\pi<\theta<\pi \} \cup \left\{ re^{i\pi} : 1<r<\infty \right\}.
\end{align}
For example, the gamma function admits the integral representation \cite[Eqs. 5.5.3 and 5.9.2]{NIST}
\begin{align}\label{}
  \Gamma(s) = \frac{1}{e^{i\pi s} - e^{-i\pi s}} \int_{H_1} e^z z^{s-1} dz, \qquad s \in \C.
\end{align}
Similarly, it is shown in \cite{FKK} that the following identity is valid:
\begin{align}\label{}
  \sum_{k=0}^\infty \frac{(-i)^k}{k!(k+ix)} = \frac{1}{e^{-\pi x} - e^{\pi x}} \int_{H_1^1} e^{iz} z^{ix -1} dz, \qquad x \in \R,
\end{align}
where $H_1^1$ denotes a contour similar to $H_1$ but the endpoints are $e^{i\pi}$ and $e^{-i\pi}$ instead of $\infty e^{i\pi}$ and $\infty e^{-i\pi}$.

In addition, the large $t$ asymptotic analysis of  several recently derived identities satisfied by  the Riemann zeta function, $\zeta(\sigma + it)$, $\sigma, t$ real \cite{F2019, Fpreprint}, requires computing the asymptotic form of certain integrals which also involve the Hankel contour; these integral depend on a real parameter, $\alpha$. In this work, a rigorous technique is introduced for analysing such integrals. Specifically, Theorem \ref{E4SDth}, given below, presents the large $\lambda$-asymptotics to all orders of the integral defined in (\ref{Ldef}): Equation (\ref{Lasymptotics}) shows that, depending on the value of $\alpha$, the relevant formula involves an asymptotic series of explicit terms, as well a term containing a specific integral.

The above situation is conceptually similar with the one occurring in the asymptotic evaluation of integrals via the steepest descent technique. In the latter case, in general, the relevant formula involves an asymptotic series of explicit terms. However, if the stationary point coincides with one of the boundaries of the integral under consideration, then the series must be supplemented with a specific integral. The latter integral is often referred to as the Bleistein integral, since the mathematician Norman Bleistein was the first to present the relevant formula \cite{BH2010}. It is shown in Lemma \ref{ABClemma} that, remarkably, the leading order of the Bleistein formula can be expressed in a form that coincides with the leading order of the integral appearing in formula (\ref{Kdef}). It is shown in \cite{FKK} that the latter formulation has certain advantages in comparison with the classical expression.

The proof of Theorem 1.1 is presented in section \ref{proofsec}. 
The leading behavior of the integral appearing in (\ref{Kdef}) is determined in section \ref{Ksec}.
The relationship between the leading form of the Bleistein formula and equation (\ref{Kdef}) is presented in section \ref{relationsec}.

Unless stated otherwise, the principal branch is used for the logarithm $\ln z$.

\begin{theorem}\label{E4SDth}
Let
\begin{align}\label{Ldef}
  L(\alpha, \lambda) = \int_{H_1} f(\alpha, w) e^{\lambda S(\alpha, w)} dw, 
\end{align}
where 
\begin{align}
 f(\alpha, w) = \frac{1}{\pi w(\frac{\pi}{2} - i\ln(w/\alpha))}, \quad
S(\alpha, w) = w - \frac{\pi}{2} + i \ln(w/\alpha).
\end{align}
Let $g(\alpha,r)$ and $g_{\reg}(\alpha,r)$ be defined by 
\begin{align}\nonumber
g(\alpha, r) & = \frac{e^{\frac{\pi i}{4}} r}{\pi(\frac{\pi}{2} - i\ln\frac{w}{\alpha})(1-iw)},
	\\ \label{gregdef}
g_{\reg}(\alpha, r) & = \frac{w e^{\frac{\pi i}{4}} r}{\pi(1-iw)}\bigg(\frac{1}{w(\frac{\pi}{2} - i\ln\frac{w}{\alpha})} - \frac{i}{w + i\alpha}\bigg),
\end{align}
with the variable $r = r(w)$ defined by
\begin{align}\label{rdef}
S(\alpha, w) - S(\alpha, -i) = -\frac{1}{2}r^2, \quad \text{i.e.} \quad
r = e^{-\frac{\pi i}{4}}\sqrt{2(1-iw+\ln(iw))}.
\end{align}
The branches in (\ref{rdef}) are such that $r = 0$ and $\frac{dr}{dw} = e^{-\frac{i\pi}{4}}$ at $w = -i$, and such that $w$ depends continuously on $r \in \R$. Define $K(\alpha, \lambda)$ by
\begin{align}\label{Kdef}
K(\alpha, \lambda) = - \frac{i e^{-\lambda i(1 + \ln \alpha)}}{\pi} \int_{-\infty}^\infty  \frac{e^{-\frac{\lambda}{2} r^2} }{w + i\alpha} \frac{wr}{w + i} dr.
\end{align}

Then, for each integer $N \geq 1$, $L$ satisfies the following two asymptotic formulas:
\begin{align}\nonumber
L(\alpha, \lambda) 
= &\; e^{\frac{\pi i}{4}} e^{-\lambda i(1 + \ln \alpha)} \sum_{j=0}^{\lfloor N/2 \rfloor} \frac{\partial_r^{2j}g_{\reg}(\alpha, 0)}{(2j)!}\Gamma \left(j+\frac{1}{2}\right)
\bigg(\frac{2}{\lambda}\bigg)^{j+\frac{1}{2}} 
	\\ \label{Lasymptotics}
& + K(\alpha, \lambda)
+ O\big(\lambda^{-\frac{N+2}{2}}\big), \qquad \lambda \to \infty, \ \ 
\alpha > 0,
\end{align}
and
\begin{align}\nonumber
L(\alpha, \lambda) = &\; e^{\frac{\pi i}{4}} e^{-\lambda i(1 + \ln \alpha)} \sum_{j=0}^{\lfloor N/2\rfloor} \frac{\partial_r^{2j}g(\alpha, 0)}{(2j)!}
\Gamma \left(j+\frac{1}{2}\right)
\bigg(\frac{2}{\lambda}\bigg)^{j+\frac{1}{2}}
	\\ \label{Lasymptoticsalphanot1}
& + O\bigg(\bigg(\frac{1}{|\alpha -1|^{N+2}} + 1\bigg) \lambda^{-\frac{N+2}{2}}\bigg), \qquad \lambda \to \infty, \ \ 
\alpha \in (0, \infty) \setminus \{1\},
\end{align}
where the error terms are uniform with respect to $\alpha$ in the given ranges, and $\lfloor N/2\rfloor$ denotes the integer part of $N/2$.

In the special case of $N = 3$, the asymptotic formulas (\ref{Lasymptotics}) and (\ref{Lasymptoticsalphanot1}) reduce to
\begin{align}\nonumber
L(\alpha, \lambda) 
= &\; e^{\frac{\pi i}{4}} e^{-\lambda i(1 + \ln \alpha)} 
\bigg\{
\frac{\alpha -1 - \ln{\alpha}}{\sqrt{\pi } (\alpha -1) \ln{\alpha} }
\bigg(\frac{2}{\lambda}\bigg)^{\frac{1}{2}} 
+ i \frac{\frac{1 + 10 \alpha  + \alpha^2}{(\alpha -1)^3}-\frac{12 - 6 (\ln \alpha) + (\ln \alpha)^2}{(\ln \alpha )^3}}{24 \sqrt{\pi }}
\bigg(\frac{2}{\lambda}\bigg)^{\frac{3}{2}} 
\bigg\}
	\\ \label{LasymptoticsN3}
& + K(\alpha, \lambda)
+ O\big(\lambda^{-\frac{5}{2}}\big), \qquad \lambda \to \infty, \ \ 
\alpha > 0,
\end{align}
and
\begin{align}\nonumber
L(\alpha, \lambda) = &\; e^{\frac{\pi i}{4}} e^{-\lambda i(1 + \ln \alpha)} 
\bigg\{
\frac{1}{\sqrt{\pi } \ln{\alpha} }
\bigg(\frac{2}{\lambda}\bigg)^{\frac{1}{2}}
-i \frac{12 - 6 (\ln \alpha) + (\ln \alpha)^2}{24 \sqrt{\pi } (\ln \alpha )^3}
\bigg(\frac{2}{\lambda}\bigg)^{\frac{3}{2}}
\bigg\}
	\\ \label{Lasymptoticsalphanot1N3}
& + O\bigg(\bigg(\frac{1}{|\alpha -1|^{5}} + 1\bigg) \lambda^{-\frac{5}{2}}\bigg), \qquad \lambda \to \infty, \ \ 
\alpha \in (0, \infty) \setminus \{1\},
\end{align}
as $\lambda \to \infty$ uniformly for $\alpha$ in the given ranges.
\end{theorem}

\section{Proof of Theorem \ref{E4SDth}}\label{proofsec}

Our goal is to find the asymptotic behavior of $L(\alpha, \lambda)$ as $\lambda \to \infty$. 

The function $S(\alpha, w)$ is analytic for $w \in \C \setminus (-\infty, 0]$, whereas $f(\alpha, w)$ is analytic for $w \in \C \setminus (-\infty, 0] \cup \{-i\alpha\}$. 
At $w = -i\alpha$, $f(\alpha, w)$ has a simple pole:
\begin{align}
f(\alpha, w) = \frac{i}{\pi(w+i\alpha)} + \frac{1}{2\pi \alpha} + O(w+i\alpha), \qquad w \to -i\alpha.
\end{align}

\begin{figure}
\begin{center}
\begin{overpic}[width=.6\textwidth]{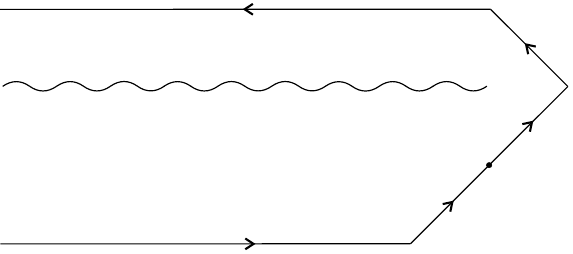}
     \put(94,37){$C_1$}
     \put(87,28){$0$}
     \put(85,12){$-i$}
\end{overpic}
     \begin{figuretext}\label{C.pdf}
         The contour $C_1$.
     \end{figuretext}
  \end{center}
\end{figure}

Noting that
$$\partial_wS(\alpha, w) = \frac{i+w}{w},$$
it follows that $L$ has a single critical point at $w = -i$. 
Let $\epsilon > 0$ be small. Let $D_\epsilon(-i)$ be the open disk of radius $\epsilon$ centered at $-i$. We deform the contour  $H_1$ to the contour $C_1$ consisting of the union of straight-line segments going from $-\infty-2i$ to $-1-2i$ to $1$ to $i$ to $-\infty + i$, see Figure \ref{C.pdf}. Then the intersection of $C_1$ with $D_\epsilon(-i)$ is the straight line segment from $-i - \epsilon e^{\frac{\pi i}{4}}$ to $-i + \epsilon e^{\frac{\pi i}{4}}$.

Letting $w = |w|e^{i\arg w}$, we find
\begin{align}
\re S(\alpha, w) = -\frac{\pi}{2} + \re w - \arg w.
\end{align}
In particular, $\re S(\alpha,w)$ is independent of $\alpha$.
Also,
\begin{align}
|f(\alpha, w)| = \frac{1}{\pi |w| \sqrt{(\frac{\pi}{2} + \arg w)^2 + (\ln\frac{|w|}{\alpha})^2}}.
\end{align}
Hence there exist constants $c, C > 0$ (depending only on $\epsilon$) such that
$$\re S(\alpha, w) \leq -c \quad \text{and} \quad |f(\alpha, w)| \leq C$$
for all $w \in C \setminus D_\epsilon(-i)$ and all $\alpha > 0$.
Thus,
\begin{align}\label{LJ}
L(\alpha, \lambda) = J(\alpha,\lambda) + O(e^{-c\lambda}), \qquad \lambda \to \infty,
\end{align}
where the error term is uniform with respect to $\alpha > 0$, and $J$ is defined by
\begin{align}\label{Jdef}
  J(\alpha,\lambda) & = \int_{C_1 \cap D_\epsilon(-i)} f(\alpha, w) e^{\lambda S(\alpha, w)} dw.
\end{align}

\subsection{Proof of (\ref{Lasymptoticsalphanot1}): the case $\alpha \in (0, \infty) \setminus \{1\}$}

Suppose $\alpha \in (0, \infty) \setminus \{1\}$ so that the pole of $f(\alpha, w)$ does not coincide with the critical point at $-i$. 
The method of steepest descent can be applied to the integral $J(\alpha, \lambda)$ as follows. 

We have
$$S(\alpha, w) = S(\alpha, -i) + \frac{i}{2}(w+i)^2 + \frac{1}{3}(w+i)^3 + O((w+i)^4), \qquad w \to -i.$$
Hence it is possible to define a new variable $z = z(w)$ by
\begin{align}\label{zdef}
S(\alpha, w) - S(\alpha, -i) = \frac{i}{2}z^2, \quad \text{i.e.} \quad
z = \sqrt{2(1-iw+\ln(iw))},
\end{align}
where the branch of the square root is fixed by the requirement that $dz/dw = 1$ at $w = -i$.
Shrinking $\epsilon > 0$ if necessary, the map $w \mapsto z$ is a conformal bijection from $D_\epsilon(-i)$ onto an open neighborhood of the origin in the complex $z$-plane. By deforming the contour $C_1$ slightly, we can assume that $C_1 \cap D_\epsilon(-i)$ is mapped onto the straight line segment from  $-\epsilon e^{\frac{\pi i}{4}}$ to $\epsilon e^{\frac{\pi i}{4}}$ in the $z$-plane.

Introducing $z$ by (\ref{zdef}) and then letting $z = r e^{\frac{\pi i}{4}}$, we find
\begin{align}\label{Jexpression}
J(\alpha, \lambda) 
= e^{\lambda S(\alpha, -i)} \int_{-\epsilon e^{\frac{\pi i}{4}}}^{\epsilon e^{\frac{\pi i}{4}}} f(\alpha, w) e^{\lambda \frac{iz^2}{2}} \frac{dw}{dz} dz
= e^{\frac{\pi i}{4}} e^{\lambda S(\alpha, -i)} \int_{-\epsilon}^\epsilon 
g(\alpha,r) e^{-\lambda \frac{r^2}{2}}  dr,
\end{align}
where
\begin{align}\label{dwdz}
\frac{dw}{dz} = \frac{iwz}{w + i}
\end{align}
and
\begin{align}\label{gdef}
g(\alpha, r) := f(\alpha, w) \frac{dw}{dz} = \frac{r e^{\frac{\pi i}{4}}}{\pi(\frac{\pi}{2} - i\ln\frac{w}{\alpha})(1-iw)}.
\end{align}
Now
$$\frac{\partial}{\partial r} = \frac{d w}{d r} \frac{\partial}{\partial w}$$
where
\begin{align}\label{dwdr}
\frac{d w}{d r} = -\frac{w r}{w + i}
\end{align}
is an analytic function of $w \in D_{2\epsilon}(-i)$. Hence, 
$$\bigg|\frac{d^j w}{d r^j}\bigg| \leq C, \qquad w \in D_\epsilon(-i), \ j = 0,1, \dots, N.$$
Also, for each $j = 0,1, \dots, N$,
\begin{align*}
 \bigg|\frac{\partial^j}{\partial w^j} \frac{1}{\frac{\pi}{2} - i\ln\frac{w}{\alpha}}\bigg|
& \leq \frac{C}{|\frac{\pi}{2} - i\ln\frac{w}{\alpha}|^{j+1}}
\leq \frac{C}{(|\re w| + ||w| - \alpha|)^{j+1}} + C
	\\
& \leq \frac{C}{(|\alpha -1| + |w + i|)^{j+1}} + C, \qquad w \in C_1 \cap D_\epsilon(-i), \ \alpha > 0.
\end{align*}
For $j = 0,1, \dots, N+1$, it follows that 
\begin{align}\label{partialgbound}
|\partial_r^j g(\alpha, r)|
\leq \frac{C}{(|\alpha -1| + |w + i|)^{j+1}} + C, \qquad
w \in C_1 \cap D_\epsilon(-i), \ \alpha > 0.
\end{align}
In particular,
\begin{align}\label{partialgboundat0}
|\partial_r^j g(\alpha, 0)|
\leq \frac{C}{|\alpha -1|^{j+1}} + C, \qquad
w \in C_1 \cap D_\epsilon(-i), \  \alpha > 0, \ j = 0,1, \dots, N.
\end{align}

Expanding $g(\alpha, r)$ in a power series at $r = 0$, we can write
\begin{align}\label{gtaylor}
g(\alpha, r) = \sum_{j=0}^N \frac{\partial_r^jg(\alpha, 0)}{j!} r^j + R_N(\alpha,r),
\end{align}
where the error term is given by
$$R_N(\alpha, r) = \int_0^r \frac{\partial_r^{N+1} g(\alpha, \tilde{r})}{N!} (r - \tilde{r})^N d\tilde{r}.$$ 

The sum in (\ref{gtaylor}) makes the following contribution to $J(\alpha, \lambda)$:
\begin{align*}
J_{\mathrm{sum}}(\alpha&, \lambda) 
= e^{\frac{\pi i}{4}} e^{\lambda S(\alpha, -i)}
\sum_{j=0}^N \frac{\partial_r^jg(\alpha, 0)}{j!}
 \int_{-\epsilon}^\epsilon  r^j e^{-\lambda \frac{r^2}{2}}  dr
	\\
& = e^{\frac{\pi i}{4}} e^{\lambda S(\alpha, -i)} \sum_{j=0}^N \frac{\partial_r^jg(\alpha, 0)}{j!}
\int_{-\infty}^\infty  r^j e^{-\lambda \frac{r^2}{2}}  dr
+ O\bigg(\bigg(\frac{1}{|\alpha - 1|^{N+1}} +1\bigg)e^{-c \lambda}\bigg)
	\\
& = e^{\frac{\pi i}{4}} e^{\lambda S(\alpha, -i)} \sum_{j=0}^{\lfloor N/2 \rfloor} \frac{\partial_r^{2j}g(\alpha, 0)}{(2j)!}
\int_{-\infty}^\infty  r^{2j} e^{-\lambda \frac{r^2}{2}}  dr
+ O\bigg(\bigg(\frac{1}{|\alpha - 1|^{N+1}} +1\bigg)e^{-c \lambda}\bigg)
	\\
& = e^{\frac{\pi i}{4}} e^{\lambda S(\alpha, -i)} \sum_{j=0}^{\lfloor N/2 \rfloor} \frac{\partial_r^{2j}g(\alpha, 0)}{(2j)!}
\bigg(\frac{2}{\lambda}\bigg)^{j+\frac{1}{2}}
\Gamma \bigg(j+\frac{1}{2}\bigg)
+ O\bigg(\bigg(\frac{1}{|\alpha - 1|^{N+1}} +1\bigg)e^{-c \lambda}\bigg)
\end{align*}
as $\lambda \to \infty$ uniformly for $\alpha \in (0, \infty) \setminus \{1\}$, where the second equality follows from (\ref{partialgboundat0}).
On the other hand, by (\ref{partialgbound}),
\begin{align*}
|R_N(\alpha, r)| 
& \leq C \int_0^r \bigg(\frac{1}{|\alpha -1|^{N+2}} + 1\bigg)(r - \tilde{r})^Nd\tilde{r}
	\\
& \leq C\bigg(\frac{1}{|\alpha -1|^{N+2}} + 1\bigg)|r|^{N+1}, \qquad \alpha > 0, \ r \in (-\epsilon, \epsilon).
\end{align*}
Hence the error term $R_N(\alpha, r)$  in (\ref{gtaylor}) makes the following contribution to $J(\alpha, \lambda)$:
\begin{align*}
J_{\mathrm{error}}(\alpha, \lambda) 
& = e^{\frac{\pi i}{4}} e^{\lambda S(\alpha, -i)} \int_{-\epsilon}^\epsilon 
R_N(\alpha,r) e^{-\lambda \frac{r^2}{2}}  dr
	\\
& = O\bigg(\int_{-\epsilon}^\epsilon 
\bigg(\frac{1}{|\alpha -1|^{N+2}} + 1\bigg)|r|^{N+1} e^{-\lambda \frac{r^2}{2}}  dr\bigg)
	\\
& = O\bigg(\bigg(\frac{1}{|\alpha -1|^{N+2}} + 1\bigg) \lambda^{-\frac{N+2}{2}}\bigg)
\end{align*}
as $\lambda \to \infty$ uniformly for $\alpha \in (0, \infty) \setminus \{1\}$.

In summary, for each integer $N \geq 1$, $J= J_{\mathrm{sum}} + J_{\mathrm{error}}$ satisfies
\begin{align}\nonumber
J(\alpha, \lambda)
= &\; e^{\frac{\pi i}{4}} e^{\lambda S(\alpha, -i)} \sum_{j=0}^{\lfloor N/2 \rfloor} \frac{\partial_r^{2j}g(\alpha, 0)}{(2j)!}
\bigg(\frac{2}{\lambda}\bigg)^{j+\frac{1}{2}}
\Gamma \left(j+\frac{1}{2}\right)
	\\ \label{Jasymptotics}
& + O\bigg(\bigg(\frac{1}{|\alpha -1|^{N+2}} + 1\bigg) \lambda^{-\frac{N+2}{2}}\bigg), \qquad \lambda \to \infty, \ \ 
\alpha \in (0, \infty) \setminus \{1\},
\end{align}
where the error term is uniform with respect to $\alpha$ in the given range. In light of (\ref{LJ}), this proves (\ref{Lasymptoticsalphanot1}).

\subsection{Proof of (\ref{Lasymptotics}): the case $\alpha \in (0, \infty)$}
 
According to (\ref{LJ}), we can rewrite $L$  in the form
\begin{align}\label{LJregI}
L(\alpha, \lambda) = J_{\reg}(\alpha, \lambda) + I(\alpha, \lambda) + O(e^{-c\lambda}), \qquad \lambda \to \infty,
\end{align}
where $J_{\reg}$ and $I$ are defined by
\begin{align}
  I(\alpha,\lambda) & = \int_{C_1 \cap D_\epsilon(-i)} \frac{i e^{\lambda S(\alpha, w)}}{\pi(w+i\alpha)} dw,
 	\\  
  J_{\reg}(\alpha,\lambda) & = \int_{C_1 \cap D_\epsilon(-i)} f_{\reg}(\alpha, w) e^{\lambda S(\alpha, w)} dw,
\end{align}
with
$$f_{\reg}(\alpha, w) := f(\alpha, w) - \frac{i}{\pi(w+i\alpha)}.$$
The regularized function $f_{\reg}$ is an analytic function of $w \in D_{2\epsilon}(-i)$ for each $\alpha > 0$. 
Hence, $J_{\reg}(\alpha, \lambda)$ can be computed to all orders via the method of steepest descent. In fact, the same arguments that led to (\ref{Jasymptotics}) with $f$ replaced with $f_{\reg}$ show that
\begin{align}\nonumber
J_{\reg}(\alpha, \lambda)
= &\; e^{\frac{\pi i}{4}} e^{\lambda S(\alpha, -i)} \sum_{j=0}^{\lfloor N/2 \rfloor} \frac{\partial_r^{2j}g_{\reg}(\alpha, 0)}{(2j)!}
\bigg(\frac{2}{\lambda}\bigg)^{j+\frac{1}{2}}
\Gamma \left(j+\frac{1}{2}\right)
	\\ \label{intfregasymptotics}
& + O\big(\lambda^{-\frac{N+2}{2}}\big), \qquad \lambda \to \infty, \ 
\alpha > 0,
\end{align}
where the error term is uniform with respect to all $\alpha > 0$ and
\begin{align}
g_{\reg}(\alpha, r) := f_{\reg}(\alpha, w) \frac{dw}{dz} 
= \frac{w e^{\frac{\pi i}{4}} r}{\pi(1-iw)}\bigg(\frac{1}{w(\frac{\pi}{2} - i\ln\frac{w}{\alpha})} - \frac{i}{w + i\alpha}\bigg).
\end{align}


Therefore it only remains to consider the asymptotics of the integral $I(\alpha, \lambda)$. 
As above, we change variables from $w$ to $r$ using (\ref{zdef}) and $z = r e^{\frac{\pi i}{4}}$. This yields
$$I(\alpha, \lambda) 
= \frac{ie^{\lambda S(\alpha, -i)}}{\pi} \int_{-\epsilon}^{\epsilon}  \frac{e^{-\lambda \frac{r^2}{2}}}{w + i\alpha} \frac{dw}{dr} dr
= -\frac{ie^{\lambda S(\alpha, -i)}}{\pi} \int_{-\epsilon}^{\epsilon}  \frac{e^{-\lambda \frac{r^2}{2}}}{w + i\alpha} \frac{w r}{w + i} dr.$$
Since $w + i\alpha$ and $w + i$ are uniformly bounded away from zero for $r \in \R \setminus [-\epsilon, \epsilon]$ and $\alpha > 0$, it follows that we can replace the contour  from  $-\epsilon$ to $\epsilon$ with the real line $\R$ with only an exponentially small error of the form  $O(e^{-c \lambda})$. Hence $I(\alpha, \lambda) 
 = K(\alpha, \lambda) + O(e^{-c \lambda})$, where $K(\alpha, \lambda)$ is given by 
\begin{align}\label{Kdef2}
K(\alpha, \lambda) = 
- \frac{ie^{\lambda S(\alpha, -i)}}{\pi}  \int_{- \infty}^{\infty}  
\frac{e^{-\lambda \frac{r^2}{2}}}{w + i\alpha} \frac{w r}{w + i}  dr.
\end{align}
Using that $S(\alpha, -i) = -i(1 + \ln \alpha)$, we see that $K$ can be expressed as in (\ref{Kdef}).
This shows that
\begin{align*}
L(\lambda, \alpha)
= &\; J_{\reg}(\alpha, \lambda)
 + K(\alpha, \lambda)
 + O(e^{-c\lambda}), \qquad \lambda \to \infty,
 \end{align*}
uniformly for $\alpha > 0$. Recalling the asymptotics of the integral $J_{\reg}(\alpha, \lambda)$ obtained in (\ref{intfregasymptotics}), this yields (\ref{Lasymptotics}).

\subsection{The special case $N = 3$}
Taylor expanding the expression for $r$ in (\ref{rdef}) around $w = -i$, and then inverting the resulting series, we obtain
\begin{align}\label{wrser}
w =  
-i + e^{\frac{i\pi}{4}} r
-\frac{1}{3} r^2
-\frac{e^{\frac{3i\pi}{4}}}{36}  r^3
+ O(r^4) \qquad \text{as $r \to 0$}.
\end{align}
Substituting this expansion into (\ref{gregdef}), we find, as $r \to 0$,
\begin{align} \nonumber
g(\alpha, r) = 
&\; \frac{1}{\pi \ln \alpha} - \frac{e^{\frac{3\pi i}{4}}(\ln \alpha -3)}{3\pi (\ln \alpha)^2}r
 - i \frac{12 - 6 (\ln \alpha) + (\ln \alpha)^2}{12 \pi (\ln \alpha)^3}r^2
+ O(r^3), 
	\\ \nonumber
g_{\reg}(\alpha, r) = &\;
\frac{\alpha -1 - \ln{\alpha}}{\pi (\alpha -1) \ln{\alpha}}
+ e^{\frac{3\pi i}{4}} 
\frac{\frac{3-\ln{\alpha}}{(\ln \alpha )^2}-\frac{2 \alpha + 1}{(\alpha -1)^2}}{3 \pi }
 + i\frac{\frac{1 + 10 \alpha  + \alpha^2}{(\alpha -1)^3}-\frac{12 - 6 (\ln \alpha) + (\ln \alpha)^2}{(\ln \alpha )^3}}{12 \pi } r^2	
	\\ \label{ggregrnear0}
& + O(r^3).
\end{align}
Using these expansions to compute the derivatives of $g$ and  $g_{\reg}$ in (\ref{Lasymptotics}) and (\ref{Lasymptoticsalphanot1}) in the special case when $N=3$, we obtain (\ref{LasymptoticsN3}) and (\ref{Lasymptoticsalphanot1N3}).

\section{Leading behavior of $K$}\label{Ksec}
According to (\ref{wrser}), we have $w = -i + e^{\frac{i\pi}{4}} r + O(r^2)$, and hence $\frac{wr}{w + i} =  -e^{\frac{i\pi}{4}} + O(r)$, as $r \to 0$. It follows that the leading order behavior of the integral $K$ in (\ref{Kdef}) is given by
\begin{align}\label{K0def}
K_0(\alpha, \lambda) & :=
\frac{ie^{-\lambda i(1 + \ln \alpha)}}{\pi} \int_{-\infty}^\infty  \frac{e^{-\frac{\lambda}{2} r^2} }{r  + (\alpha - 1)e^{\frac{\pi i}{4}}} dr.
\end{align}
We can estimate the error in the approximation $K \approx K_0$ by noting that
\begin{align}\nonumber
\frac{1}{w + i\alpha} 
& = \frac{1}{-i + e^{\frac{i\pi}{4}} r + O(r^2) + i\alpha} 
= \frac{1}{e^{\frac{i\pi}{4}} r + i(\alpha -1)} \frac{1}{1 + O(\frac{r^2}{e^{\frac{i\pi}{4}} r + i(\alpha -1)})}
	\\ \label{1wialphaestimate}
& = \frac{1}{e^{\frac{i\pi}{4}} r + i(\alpha -1)}(1 + O(r))
= \frac{1}{e^{\frac{i\pi}{4}} r + i(\alpha -1)} + O\Big(\frac{r}{\alpha -1}\Big)
\end{align}
as $r \to 0$ uniformly for $\alpha \in (0,+\infty)\setminus\{1\}$.
By (\ref{Kdef}) and (\ref{K0def}), we have
\begin{align*}
K(\alpha, \lambda) -  K_0(\alpha, \lambda)
= &
-\frac{i e^{-\lambda i(1 + \ln \alpha)}}{\pi} \int_{-\infty}^\infty e^{-\frac{\lambda}{2} r^2} \frac{1}{w + i\alpha} \bigg(\frac{wr}{w + i} + e^{\frac{i\pi}{4}}\bigg) dr
	\\
& + \frac{i e^{-\lambda i(1 + \ln \alpha)}}{\pi} \int_{-\infty}^\infty e^{-\frac{\lambda}{2} r^2} \bigg(\frac{1}{w + i\alpha} - \frac{1}{e^{\frac{i\pi}{4}} r + i(\alpha -1)}\bigg) e^{\frac{i\pi}{4}} dr.
\end{align*}
Hence, if $\epsilon > 0$ is sufficiently small, (\ref{1wialphaestimate}) implies that
\begin{align*}
|K(\alpha, \lambda) -  K_0(\alpha, \lambda) |
\leq &\; Ce^{-c \lambda \epsilon^2}
+ C \int_{-\epsilon}^\epsilon e^{-\frac{\lambda}{2} r^2} \Big|\frac{1}{w + i\alpha }\Big|  |r| dr
+ C \int_{-\epsilon}^\epsilon e^{-\frac{\lambda}{2} r^2} \Big|\frac{r}{\alpha -1}\Big|dr
	\\
\leq&\; Ce^{-c \lambda \epsilon^2} + C \int_{-\epsilon}^\epsilon e^{-\frac{\lambda}{2} r^2} \Big|\frac{r}{\alpha -1}\Big|dr
\leq \frac{C}{\lambda |\alpha -1|} 
\end{align*}
for all sufficiently large $\lambda > 0$ and all $\alpha \in (0,+\infty)\setminus\{1\}$.

\section{Relation to the Bleistein formula}\label{relationsec}

It is well known that if $\zeta^*$ is a stationary point, then as $\zeta^*$ approaches the boundary point $p$, the leading contribution as $t \to \infty$ of the integral 
$$\int_p^{\infty e^{i\phi}} e^{-it g(\zeta)} d\zeta, \qquad t\in \R,$$
which will be denoted by $J_0$, is given by the Bleistein formula
\begin{align}\label{bleistein}
J_0 = \frac{2A}{g_\zeta(p)} e^{-itg(\zeta^*)} \frac{B(\sqrt{t} A)}{\sqrt{t}},
\end{align}
with $A$ and $B$ defined by
\begin{align}\label{ABdef}
A = \sqrt{g(p) - g(\zeta^*)}, \qquad B(z) = \int_z^{e^{-\frac{i\pi}{4}}\infty} e^{-i\zeta^2} d\zeta.
\end{align}
Remarkably, $B$ can be rewritten in terms of $K_0$, where $K_0$ denotes the leading asymptotics of $K$ in (\ref{K0def}).
Indeed, performing the change of variables $\zeta = \sqrt{\lambda} r$ in (\ref{K0def}), we observe that
$$ K_0(\alpha, \lambda) = - \frac{2 e^{\frac{i\pi}{4}}}{\sqrt{\pi}} e^{-\lambda i(1 + \ln \alpha)} \hat{K}_0\bigg(i(\alpha - 1)\sqrt{\frac{\lambda}{2}} \bigg),$$
where 
$$\hat{K}_0(z) := - \frac{e^{\frac{i\pi}{4}} }{2\sqrt{\pi}} \int_{-\infty}^\infty  \frac{e^{-\frac{1}{2} \zeta^2} }{\zeta  + \sqrt{2} e^{-\frac{\pi i}{4}} z} d\zeta.$$
The next lemma shows that the function $B(z)$ appearing in the Bleistein formula (\ref{bleistein}) can be expressed in terms of $\hat{K}_0(z)$.

\begin{lemma}\label{ABClemma}
The function $B(z)$  defined in (\ref{ABdef}) satisfies
\begin{align}\label{BofA}
B(z) = \begin{cases} e^{-iz^2}\hat{K}_0(z), & \im(e^{-\frac{i\pi}{4}} z) < 0, 
	\\
e^{-iz^2} \hat{K}_0(z) + \sqrt{\pi} e^{-\frac{i\pi}{4}} , & \im(e^{-\frac{i\pi}{4}} z) > 0.
\end{cases}
\end{align}
\end{lemma}
\proofbegin
Define $\varphi(z)$ for $z \in \C \setminus (e^{\frac{i\pi}{4}}\R)$ by $\varphi(z) = e^{-iz^2} \hat{K}_0(z)$.
For $z \in \C \setminus (e^{\frac{i\pi}{4}}\R)$, $\varphi(z)$ is analytic and satisfies
\begin{align}
  \varphi'(z) = -2iz \varphi(z) + \frac{e^{-iz^2}}{\sqrt{2\pi}} \tilde{\varphi}(z),
\quad \text{where}\quad
  \tilde{\varphi}(z) := \int_{-\infty}^{\infty} \frac{e^{-\frac{\zeta^2}{2}} d\zeta}{(\zeta + \sqrt{2} e^{-\frac{i\pi}{4}} z)^2}.
\end{align}
Using integration by parts it follows that, for $e^{-\frac{i\pi}{4}} z \in \C \setminus \R$,
\begin{align}\nonumber
  \tilde{\varphi}(z)
  & = -\int_{-\infty}^{\infty} \frac{e^{-\frac{\zeta^2}{2}} \zeta d\zeta}{\zeta + \sqrt{2} e^{-\frac{i\pi}{4}} z}
  =  -\int_{-\infty}^{\infty} e^{-\frac{\zeta^2}{2}} d\zeta 
  + \sqrt{2} e^{-\frac{i\pi}{4}} z \int_{-\infty}^{\infty} \frac{e^{-\frac{\zeta^2}{2}} d\zeta}{\zeta + \sqrt{2} e^{-\frac{i\pi}{4}} z}
	\\ \label{tildeCsimplified}
 & =  - \sqrt{2\pi}
  - \sqrt{2} e^{-\frac{i\pi}{4}} z  \frac{2\sqrt{\pi}}{e^{\frac{i\pi}{4}} e^{-iz^2}} \varphi(z).
\end{align}
Hence
$$\varphi'(z) = -2iz \varphi(z) + \frac{e^{-iz^2}}{\sqrt{2\pi}} \Big( - \sqrt{2\pi} - \sqrt{2} e^{-\frac{i\pi}{4}} z  \frac{2\sqrt{\pi}}{e^{\frac{i\pi}{4}} e^{-iz^2}} \varphi(z)\Big) 
 = - e^{-iz^2}.$$
which shows that $B'(z) = \varphi'(z)$ for $e^{-\frac{i\pi}{4}} z \in \C \setminus \R$. It follows that 
\begin{align}
B(z) = \varphi(z) + \begin{cases} c_1, & \im(e^{-\frac{i\pi}{4}} z) < 0, 
	\\
c_2, & \im(e^{-\frac{i\pi}{4}} z) > 0.
\end{cases}
\end{align}
Since both $B(z)$ and $\varphi(z)$ tend to $0$ as $z$ tends to $e^{-\frac{i\pi}{4}} \infty$, we see that $c_1 = 0$.
On the other hand, for $e^{-\frac{i\pi}{4}} z \in \R$, we have
$$\varphi_+(z) - \varphi_-(z) 
= 2\pi i \underset{\zeta = -\sqrt{2} e^{-\frac{i\pi}{4}} z}{\res}\bigg(\frac{e^{\frac{i\pi}{4}} e^{-iz^2}}{2\sqrt{\pi}} \frac{e^{-\frac{\zeta^2}{2}}}{\zeta + \sqrt{2} e^{-\frac{i\pi}{4}} z}\bigg)
= -\sqrt{\pi} e^{-\frac{i\pi}{4}} $$
where $\varphi_\pm(z) =  \lim_{\epsilon \to 0^+} \varphi(z \pm e^{\frac{3i\pi}{4}}\epsilon)$.
Since $B$ is continuous across the line $e^{\frac{i\pi}{4}} \R$ and $c_1=0$, we obtain $c_2 = B_+(z) - \varphi_+(z) = \varphi_-(z) - \varphi_+(z) = \sqrt{\pi} e^{-\frac{i\pi}{4}}$.
\proofend

\begin{remark}\label{Azeroremark}
For $z=0$, the function $B$ defined in (\ref{ABdef}) equals $B(0) = \sqrt{\pi} \exp(-i\pi/4)/2$.
This is consistent with the right-hand side of (\ref{BofA}) because, as $z \to 0$ with $\im(e^{-\frac{i\pi}{4}} z) \lessgtr 0$, the Plemelj formula shows that
$$\hat{K}_0(z)  \to -\frac{e^{\frac{i\pi}{4}}}{2\sqrt{\pi}} \bigg(\dashint_{-\infty}^{\infty} \frac{e^{-\frac{\zeta^2}{2}} d\zeta}{\zeta} \pm i\pi \underset{\zeta = 0}{\res} \frac{e^{-\frac{\zeta^2}{2}}}{\zeta} \bigg)
= -\frac{e^{\frac{i\pi}{4}}}{2\sqrt{\pi}}(0 \pm i\pi) = \pm \frac{\sqrt{\pi}}{2} e^{-\frac{i\pi}{4}}.$$
\end{remark}

\begin{remark}
The function $B$ can be expressed in terms of the error function:
$$B(z) = \int_z^{e^{-\frac{i\pi}{4}}\infty} e^{-i\zeta^2} d\zeta 
= \frac{1}{2} e^{\frac{3i\pi}{4}} \sqrt{\pi } \left(\erf(e^{\frac{i\pi}{4}} z)-1\right).$$
\end{remark}

\appendix
\section{Numerical verification}
\begin{figure}
\begin{center}
\begin{overpic}[width=.3\textwidth]{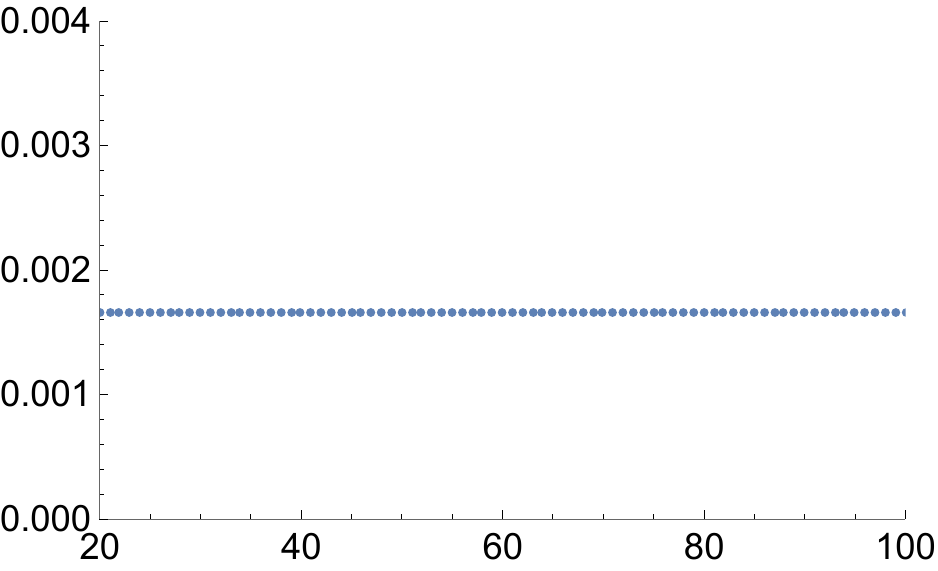}
     \put(101,4){\footnotesize $\lambda$}
     \put(1,64){\footnotesize $\lambda^{\frac{N+2}{2}} D_N$}
\end{overpic}
\hspace{4mm}
\begin{overpic}[width=.3\textwidth]{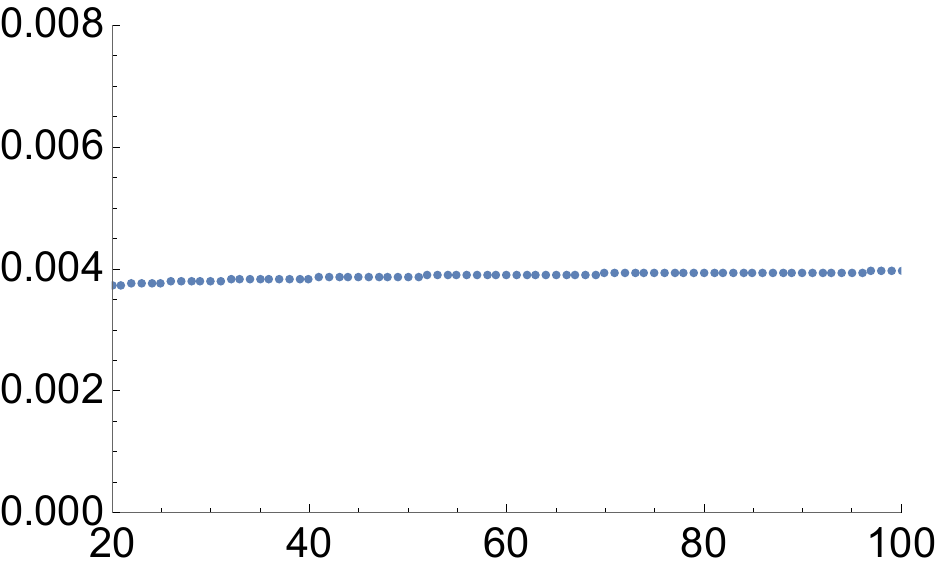}
     \put(101,4){\footnotesize $\lambda$}
     \put(1,64){\footnotesize $\lambda^{\frac{N+2}{2}} D_N$}
\end{overpic}
\hspace{4mm}
\begin{overpic}[width=.3\textwidth]{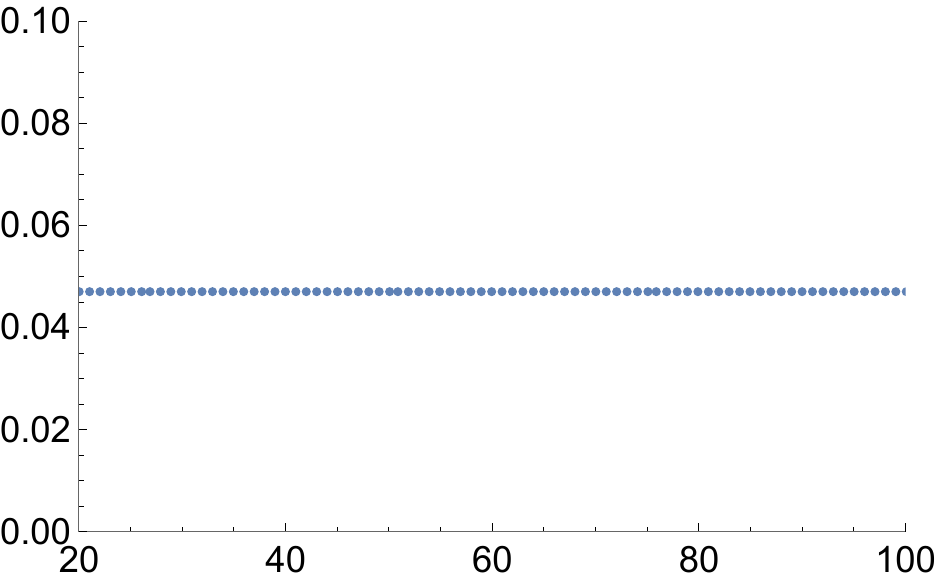}
     \put(101,4){\footnotesize $\lambda$}
     \put(1,64){\footnotesize $\lambda^{\frac{N+2}{2}} E_N$}
\end{overpic}
     \begin{figuretext}\label{numericalfig}
         Left: The error term in (\ref{Lasymptotics}) multiplied by $\lambda^{\frac{N+2}{2}}$ as a function of $\lambda \in [20,100]$ for $N = 3$ and $\alpha = 5$.
         Middle: The error term in (\ref{Lasymptotics}) multiplied by $\lambda^{\frac{N+2}{2}}$ as a function of $\lambda \in [20,100]$ for $N = 3$ and $\alpha = 1 + \frac{1}{\sqrt{\lambda}}$.
         Right: The error term in (\ref{Lasymptoticsalphanot1}) multiplied by $\lambda^{\frac{N+2}{2}}$ as a function of $\lambda \in [20,100]$ for $N = 3$ and $\alpha = 5$.
The plots are consistent with Theorem \ref{E4SDth} according to which all three quantities should be $O(1)$ as $\lambda \to \infty$.
     \end{figuretext}
  \end{center}
\end{figure}

To verify the formulas of Theorem \ref{E4SDth} numerically, we consider the differences between the left- and right-hand sides of (\ref{Lasymptotics}) and (\ref{Lasymptoticsalphanot1}), respectively:
$$D_N(\alpha, \lambda) := 
L(\alpha, \lambda) -  e^{\frac{\pi i}{4}} e^{-\lambda i(1 + \ln \alpha)} \sum_{j=0}^{\lfloor N/2 \rfloor} \frac{\partial_r^{2j}g_{\reg}(\alpha, 0)}{(2j)!}\Gamma \left(j+\frac{1}{2}\right)
\bigg(\frac{2}{\lambda}\bigg)^{j+\frac{1}{2}} 
- K(\alpha, \lambda)$$
and
$$E_N(\alpha, \lambda) := L(\alpha, \lambda) - e^{\frac{\pi i}{4}} e^{-\lambda i(1 + \ln \alpha)} \sum_{j=0}^{\lfloor N/2\rfloor} \frac{\partial_r^{2j}g(\alpha, 0)}{(2j)!}
\Gamma \left(j+\frac{1}{2}\right)
\bigg(\frac{2}{\lambda}\bigg)^{j+\frac{1}{2}}.$$
According to Theorem \ref{E4SDth}, $\lambda^{\frac{N+2}{2}} D_N(\alpha, \lambda)$ should be of order $O(1)$ as $\lambda \to \infty$, uniformly with respect to $\alpha > 0$. In Figure \ref{numericalfig}, the quantity $\lambda^{\frac{N+2}{2}} D_N(\alpha, \lambda)$ is shown as a function of $\lambda$ for a fixed value of $\alpha$ (Figure \ref{numericalfig}, left) and for a case where $\alpha$ approaches $1$ as $\lambda \to \infty$ (Figure \ref{numericalfig}, middle); in both cases, the quantity $\lambda^{\frac{N+2}{2}} D_N(\alpha, \lambda)$ appears to tend a constant for large $\lambda$ in agreement with Theorem \ref{E4SDth}.

Similarly, according to Theorem \ref{E4SDth}, $\lambda^{\frac{N+2}{2}} E_N(\alpha, \lambda)$ should be of order $O(1)$ as $\lambda \to \infty$ for any fixed $\alpha \in (0,1) \cup (1, +\infty)$. In Figure \ref{numericalfig} (right), the quantity $\lambda^{\frac{N+2}{2}} E_N(\alpha, \lambda)$ is shown as a function of $\lambda$ for a fixed value of $\alpha$; it seems to tend a constant for large $\lambda$ in agreement with Theorem \ref{E4SDth}.

 \bigskip
\noindent
{\bf Acknowledgement} {\it JL acknowledges support from the Swedish Research Council, Grant No. 2021-03877.}

\printbibliography

@article{F2019,
    author = {A. S. Fokas},
    title = {A novel approach to the Lindelöf hypothesis},
    journal = {Trans. Math. Appl.},
    volume = {3},
%    number = {1},
%    pages = {tnz006},
    year = {2019},
 %   month = {09},
 %   issn = {2398-4945},
 %   doi = {10.1093/imatrm/tnz006},
%    url = {https://doi.org/10.1093/imatrm/tnz006},
%    eprint = {https://academic.oup.com/imatrm/article-pdf/3/1/tnz006/32856499/tnz006.pdf},
}

@article{NIST,
    author = {F. W. J. Olver and A. B. Olde Daalhuis and D. W. Lozier and B. I. Schneider and R. F. Boisvert and C. W. Clark and B. R. Miller and B. V. Saunders and H. S. Cohl and M. A. McClain},
    title = {NIST Digital Library of Mathematical Functions},
    year = {Release 1.2.4 of 2025-03-15},
}

@article{Fpreprint,
    author = {A. S. Fokas},
    title = {Remarkable identities satisfied by the Riemann zeta function and explicit large-$t$ asymptotic formulas (preprint)},
}

@article{FKK,
    author = {A. S. Fokas and D. Kyriakopoulou and K. Kalimeris},
    title = {Asymptotic techniques useful for the large $t$ evaluation of the Riemann zeta function (in preparation)},
}

@Book{BH2010,
  author     = {N. Bleistein and R. A. Handelsman},
  publisher  = {Dover Publications},
  title      = {Asymptotic Expansions of Integrals},
  year       = {2010},
}

\end{document}